\documentclass{amsart}

\usepackage{latexsym}
\usepackage[dvips]{graphics}
\usepackage{psfrag}

\newtheorem{theorem}{Theorem}[section]
\newtheorem{proposition}[theorem]{Proposition}

\newtheorem{corollary}[theorem]{Corollary}
\newtheorem{lemma}[theorem]{Lemma}

\theoremstyle{definition}
\newtheorem{example}[theorem]{Example}
\newtheorem{defn}[theorem]{Definition}
\newtheorem{remark}[theorem]{Remark}

\newcommand{\es}{\emptyset}

\begin{document}

\title{Boolean Term Orders and the Root System~$B_n$}
\author{Diane Maclagan }
\email{maclagan@math.berkeley.edu} 

\address{
Diane Maclagan\\
Department of Mathematics,\\
University of California, Berkeley,\\
Berkeley, CA 94720}

\begin{abstract}
A boolean term order is a total order on subsets of $[n]=\{1,\ldots,n\}$ such that $\es \prec \alpha$ for all $\alpha \subseteq [n], \alpha \neq \es$, and $\alpha \prec \beta \Rightarrow \alpha \cup \gamma \prec \beta \cup \gamma$ for all $\gamma$ with $\gamma \cap (\alpha \cup \beta) = \es$.  Boolean term orders arise in several different areas of mathematics, including Gr\"obner basis theory for the exterior algebra, and comparative probability.  

The main result of this paper is that boolean term orders correspond
to one-element extensions of the oriented matroid $\mathcal M (B_n)$,
where $B_n$ is the root system $\{e_i:1 \leq i \leq n \} \cup \{e_i
\pm e_j :1 \leq i < j \leq n \}$.  This establishes boolean term
orders in the framework of the Baues problem, in the sense of
\cite{reiner}.  We also define a notion of coherence for a boolean term order, and a flip relation between different term orders.  Other results include examples of noncoherent term orders, including an example exhibiting flip deficiency, and enumeration of boolean term orders for small values of $n$.
\end{abstract}

\maketitle

\section{Introduction}

\begin{defn} \label{bto} A {\em boolean term order} is a total order on subsets of $[n]=\{1,\ldots,n\}$ such that:

\begin{enumerate}
\item $\es \prec \alpha$ for all $\alpha \subseteq [n], \alpha \neq
\es$.
\item If $\alpha \prec \beta$, and $\gamma \cap (
\alpha \cup \beta)=\emptyset$ then  $\alpha \cup \gamma \prec
\beta \cup \gamma$ \label{axiom2}
\end{enumerate}
Where there can be no confusion, we will use the phrase {\em term order}.
\end{defn}

Boolean term orders arise in several different areas of mathematics.
The motivating example for this paper is term orders in the exterior
algebra over a vector space of dimension $n$.  A monomial $x^{\alpha}$
in the exterior algebra corresponds to the subset of $[n]$ given by
$\mathrm{supp}(\alpha)=\{i:x_i | x^{\alpha}\}$.  The axioms above then
correspond to the standard ones for Gr\"obner basis theory term orders
(see \cite{stokes}).

Another incarnation of boolean term orders is as antisymmetric
comparative probability relations.  The basic idea there is that we
have a finite set $\Omega$ of events, and are interested in which
events are more probable than others, as opposed to the exact
probability of each event.  We associate with a subset of $\Omega$ the
(comparative) probability that at least one of the events in the
subset occurs.  If we demand that for any two subsets $A$ and $B$ of
$\Omega$ either $A \prec B$ ($B$ is more probable than $A$), or $B
\prec A$, then we have a boolean term order with $n=|\Omega|$.  An overview of the theory of comparative probability can be be found in \cite{fishburn}.

The main result of this paper is that boolean term orders correspond
to one-element extensions of the oriented matroid $\mathcal M (B_n)$,
where $B_n$ is the root system $\{e_i:1 \leq i \leq n \} \cup \{e_i
\pm e_j :1 \leq i < j \leq n \}$.  This establishes boolean term
orders in the frame work of the Baues problem, in the sense of
\cite{reiner}.  Section 2 contains precise definitions of term orders,
including a notion of coherence.  The connection to the Baues problem
is strengthened in Section 3 with the introduction of a notion of a
flip relation between different boolean term orders.  Examples are
given of term orders exhibiting flip deficiency.  Section 4 examines
the structure of the set of coherent term orders, and part of the connection
to $B_n$.  This connection is fully established in Section 5 with the
proof of the main theorem described above.  Finally, Section 6
contains some specific examples of boolean term orders, and
enumeration for small values of $n$.

\section{Definitions}

\begin{defn} \label{gpto} A {\em generalized partial term order}
 is an irreflexive partial order, $\prec$, on the set of subsets of
 $[n]$ such that:

\begin{enumerate}

\item $\alpha \prec \beta \Leftrightarrow \alpha \cup \gamma \prec
\beta \cup \gamma$  for all $\gamma$ with $\gamma \cap (\alpha \cup
\beta)=\emptyset$

\item If $\alpha$ and $\beta$ are not comparable (written
 $\alpha \sim \beta$), then $\{ \gamma : \gamma \prec \alpha \}=\{\gamma
 : \gamma \prec \beta \}$, and $\{ \gamma : \alpha \prec \gamma \}=\{\gamma :
 \beta \prec \gamma \}$
\end{enumerate}

  We write $\alpha \preceq \beta$ if $\alpha \sim \beta$ or $\alpha
 \prec \beta$.  A {\em partial boolean term order} additionally
 satisfies $\es \prec \alpha$ for all $\alpha \neq \es$.  A {\em
 generalized term order} is a generalized partial term order where
 $\alpha \sim \beta \Rightarrow \alpha=\beta$.  Note that a  boolean term
 order is a generalized term order with the additional requirement
 that $\emptyset \prec \alpha$ for all $\alpha \neq \emptyset$.
\end{defn}

\begin{lemma} \label{multlemma}
\begin{enumerate}

\item The first condition of Definition \ref{gpto} is equivalent to
requiring that whenever $a \prec b$, $c \prec d$ with $a \cap c =
\emptyset$ and $b \cap d = \emptyset$, then $a \cup c \prec b \cup d$
\label{multlemmaa}

\item Assuming the first condition of Definition \ref{gpto}, the second condition is equivalent to requiring that if either
$a \cup c \sim b \cup d$ or $a \cup c =b \cup d$, with $a \cap c = \emptyset$  and $b \cap d = \emptyset$, and $b \prec a$, then $c
\prec d$. \label{multlemmab}
\end{enumerate}
\end{lemma}
\begin{proof}
\begin{enumerate}
\item Assume $\prec$ satisfies the first condition of Definition \ref{gpto} and write
$b=b^{\prime} \cup l$, $c=c^{\prime} \cup l$, where $b^{\prime}$, 
$c^{\prime}$, and $l$ are all pairwise disjoint.   Then $a \cup c^{\prime} \prec
b \cup c^{\prime}=b^{\prime} \cup c \prec b^{\prime} \cup d$, so
$a \cup c^{\prime} \cup l \prec b^{\prime} \cup d \cup l$.  This proves one
implication.

  Conversely, assume that $\prec$ is a total order on subsets of $[n]$
satisfying the condition of the lemma, but not satisfying the first
condition of Definition \ref{gpto}, so there is a pair $a \prec b$
with $b \cup c \prec a \cup c$ for some $c$ with $c \cap (a \cup b)=
\emptyset$.  We may assume that $a \cap b = \es$, so the condition of
the lemma implies $a \cup b \cup c \prec a \cup b \cup c$, a
contradiction.

\item
Suppose $\prec$ satisfies both conditions of Definition \ref{gpto}, $b
\prec a$, and either $a \cup c \sim b \cup d$ or $a \cup c=b \cup d$.
If $d \prec c$ we would have $b \cup d \prec a \cup c$ from above. If
$c \sim d$, then write $b=b^{\prime} \cup l$, $c=c^{\prime} \cup l$ as
above.  Then $b^{\prime} \cup c=b \cup c^{\prime} \prec a \cup
c^{\prime}$.  But $b^{\prime} \cup c \sim b^{\prime} \cup d$, so
$b^{\prime} \cup d \prec a \cup c^{\prime}$.  From this contradiction
we conclude $c \prec d$.

Conversely, let $\prec$ be a partial order on monomials satisfying the
first condition of Definition \ref{gpto}, and such that whenever $b
\prec a$ and either $a\cup c \sim b \cup d$ or $a \cup c=b \cup d$, with $a \cap c = \emptyset$ and $b \cap d = \emptyset$, 
then $c \prec d$.  Suppose $\alpha \sim \beta$, and $\alpha \prec
\gamma$.  We can assume that $\alpha \cap \beta \cap \gamma=
\emptyset$.  Write
\begin{eqnarray*}
\alpha & = & \alpha^{\prime} \cup \delta  \cup  \phi \\
\beta & = & \beta^{\prime} \cup \delta \cup \psi \\
\gamma & = & \gamma^{\prime} \cup \phi \cup \psi \\
\end{eqnarray*}
where $\alpha^{\prime},\beta^{\prime}, \gamma^{\prime}, \delta, \phi,
\psi$ are all disjoint.  Then $\alpha \sim \beta$ implies $\alpha \cup
\gamma^{\prime} \sim \beta \cup \gamma^{\prime}$.  Writing this out in
full, we have $(\gamma^{\prime} \cup \phi) \cup (\alpha^{\prime} \cup
\delta) \sim (\gamma^{\prime} \cup \psi) \cup (\beta^{\prime} \cup
\delta)$.  But $\alpha \prec \gamma$ means $\alpha^{\prime} \cup \delta
\prec \gamma^{\prime} \cup \psi$, so $\beta^{\prime} \cup \delta \prec
\gamma^{\prime} \cup \phi$, and thus, multiplying by $\psi$, we get
$\beta \prec \gamma$, as required by the definition.  The case with
$\gamma \prec \alpha$ is the same with the inequalities reversed.
\end{enumerate}
\end{proof}

\begin{corollary} \label{complement}
If $\prec$ is a generalized partial term order with
$a \prec b$ then $[n] \setminus b \prec [n] \setminus a$.
\end{corollary} 

\begin{proof} We have $[n]=a \cup ([n] \setminus a)$, and $[n]=b \cup ([n]
\setminus b)$, so the result follows from part 2 of Lemma
\ref{multlemma}. \hfill $\Box$
\end{proof}

This means that the second half of the term order is the complement of
the first, in reverse order.

\begin{defn} A boolean term order is {\em coherent} if there is a weight 
vector $w=(w_1,w_2,\ldots,w_n) \in \mathbb{N}^n$ such that 
$$\alpha \prec \beta \Leftrightarrow \sum_{i \in \alpha} w_i < \sum_{j
\in \beta} w_j$$ In the interpretation of boolean term orders as
Gr\"obner basis term orders in the exterior algebra, a term order is
coherent if it can be extended to a Gr\"obner basis term order on all
the monomials in $n$ variables.  This follows from the fact that a
Gr\"obner basis term order in $n$ variables can be induced up to a
given finite degree by an integral weight vector (see, for example,
Chapter 15 of \cite{eisenbud}). In the comparative probability
language these are the comparative probability orders which have an
agreeing probability measure, and are known as additive antisymmetric
comparative probability orders.  The question of whether all
antisymmetric comparative probability orders were additive was first
raised by de Finetti in 1951 \cite{deFinetti}, and first answered in
1959 by Kraft, Pratt, and Seidenberg \cite{KPS}.
\end{defn}

\begin{example}  \label{noncohorder} A {\em noncoherent} boolean term  order for $n=5$ is:
\begin{eqnarray*} 
& &\es \prec \{1\} \prec \{2\} \prec \{3\} \prec \{4\} \prec \{1, 2\}
\prec \{5 \} \prec \{1, 3\} \prec \{2, 3\} \prec \\ & & \{1, 4\} \prec
\{1, 5\} \prec \{2, 4\} \prec \{2, 5\} \prec \{3, 4\} \prec \{1, 2,
3\} \prec \{1, 2, 4\} \prec \\ & & \{3, 5\}  \prec \{4, 5\} \prec \{1,
2, 5\} \prec \{1, 3, 4 \} \prec \{1, 3, 5\} \prec \{2, 3, 4\} \prec
\\ & & \{2, 3, 5\} \prec  \{1, 4, 5\} \prec \{2, 4, 5\} \prec \{1, 2,
3, 4\} \prec \{3, 4, 5 \} \prec \\ & &  \{1, 2, 3, 5 \} \prec  \{1, 2, 4, 5\}
\prec \{1, 3, 4, 5\} \prec \{2, 3, 4, 5\} \prec \{1, 2, 3, 4, 5 \}
\end{eqnarray*}
\end{example}

To see that this term order is not coherent, we notice that $\{4 \} \prec \{1,
2\}$, $\{2, 3 \} \prec \{ 1, 4\}$, $\{1, 5\} \prec \{2, 4\}$, and
$\{1, 2, 4\} \prec \{3, 5\}$.  In the exterior algebra
characterization, this is $x_4 \prec x_1x_2$, $x_2x_3 \prec x_ 1x_ 4$,
$x_1x_5 \prec x_2x_4$, and $x_1x_2x_4 \prec x_3x_5$.
Note that if $\prec$ extended to an order on all
monomials, we could multiply all the left and right sides to get
$x_1^2x_2^2x_3x_4^2x_5 \prec x_1^2x_2^2x_3x_4^2x_5$, a contradiction.
So there is no order on the polynomial ring extending $\prec$.

This method of giving a certificate for the noncoherency of a boolean
term order has received attention recently in the work of Fishburn
\cite{fish1}, \cite{fish2}.  An open problem is to give a sharp upper
bound on the number of inequalities needed in such a certificate for 
term orders on subsets of $[n]$.

\section{Flips for Term Orders} \label{flipsection}

In this section we place boolean term orders in the framework of the
Baues problem.  Specifically, we define a notion of flip for a term
order, analogous to bistellar flips for triangulations.

\begin{defn}  A {\em primitive pair} in a term order is a  pair $\alpha \prec
 \beta$, with $\alpha \cap \beta = \es$, which is {\em adjacent}, in
the sense that there is no $\gamma$ with $\alpha \prec \gamma \prec
\beta$.
\end{defn}

\begin{example} In the term order of Example \ref{noncohorder}, the primitive pairs are: $\es \prec \{1 \}$, $\{ 1 \} \prec
\{2 \}$, $\{2 \} \prec \{3\}$, $\{3 \} \prec \{4\}$, $\{4 \} \prec
\{1,2\}$, $\{1,2 \} \prec \{5\}$, $\{5 \} \prec \{1,3\}$, $\{2,3\}
\prec \{1,4\}$, $\{1,5\} \prec \{2,4\}$, $\{2,5\} \prec \{3,4\}$, 
and $\{1,2,4\} \prec \{3,5\}$.
\end{example}

\begin{proposition} \label{primpairs} The order on the primitive pairs
 determines the boolean term order, in the sense that if $\prec_1$ and
 $\prec_2$ are different term orders, then there is a primitive pair
of $\prec_1$, $\alpha \prec_1 \beta$, such that $\beta \prec_2
 \alpha$.
\end{proposition}

 \begin{proof} Suppose that $\prec_1$ and $\prec_2$ are two distinct term
orders which have the same order on all the primitive pairs of
$\prec_1$.  Let $\alpha$ be such that, for all $\mu, \nu \subseteq
[n]$, $\mu \prec_1 \nu \prec_1 \alpha \Leftrightarrow \mu \prec_2 \nu
\prec_2 \alpha$, and suppose that $\alpha$ is the greatest subset with respect to
$\prec_1$ with this property.  Denote the next subset for $\prec_1$ by
$\beta$, and for $\prec_2$ by $\gamma$.  By assumption $\beta \neq
\gamma$.
	
Let $\delta$ be the subset immediately preceding $\gamma$ in
$\prec_1$.  We know that $\delta \cap \gamma \neq \emptyset$, as
otherwise $\delta \prec_1 \gamma$ is a primitive pair of $\prec_1$, in which case 
from consideration of $\prec_2$ we see that $\delta$ is at most
$\alpha$, and so $\beta=\gamma$.  Denote $\delta \setminus \gamma$ and
$\gamma \setminus \delta$ by $\delta^{\prime}$ and $\gamma^{\prime}$
respectively.  Looking at $\prec_2$, we see $\gamma^{\prime}$ is at
most $\alpha$, so $\delta^{\prime} \prec_2 \gamma^{\prime}$, and thus
$\delta \prec_2 \gamma$, contradicting the adjacency of $\alpha$ and
$\gamma$ with respect to $\prec_2$.  \hfill $\Box$
\end{proof}

\begin{defn} A {\em flippable pair} in a term order $\prec$ is a primitive pair
$\alpha \prec \beta$ such that all multiples $\alpha  \cup l
\prec \beta  \cup l$ with $l \cap (\alpha \cup \beta)=\emptyset$ are adjacent.
\end{defn}

\begin{example} In the term order of
example \ref{noncohorder}, the flippable pairs are $\{4 \} \prec \{1,2\}$,
$\{2,3 \} \prec \{1,4\}$, $ \{1,5 \} \prec \{2,4\}$, $\{2,5 \} \prec \{3,4\}$,
and $\{1,2,4 \} \prec \{3,5\}$.
\end{example}

\begin{defn}  Given a boolean term order $\prec_1$ and a flippable pair $\alpha
\prec_1 \beta$, with $\alpha \neq \es$, we construct a new total order
$\prec_2$ by exchanging the order of $\alpha \cup l$ and $\beta  \cup l$
for all $l \subseteq [n] \setminus (\alpha \cup \beta)$.
We say that $\prec_2$ is obtained from $\prec_1$ by {\em flipping
across}  $\alpha \prec_1 \beta$.
\end{defn}

\begin{proposition}
The $\prec_2$ constructed above is a boolean term order.
\end{proposition}

\begin{proof} Since $\es$ is still the smallest element in $\prec_2$, we need
only check that the order satisfies the second condition of Definition
\ref{bto}.  Suppose it does not, so there exist $a, \ b, \ c \subseteq
[n]$ with $c \cap (a \cup b)=\emptyset$ such that $a \prec_2 b$, but
$b \cup c \prec_2 a \cup c$.  The only pairs whose orders have changed
from $\prec_1$ are the multiples of $\alpha \prec_1 \beta$, so one of
$a \prec_2 b$, and $b \cup c \prec_2 a \cup c$ must be such a
multiple.  But if one pair is a multiple of $\beta \prec_2 \alpha$, so
is the other, and so their order is still consistent, as it is the
reverse of the consistent order under $\prec_1$.  Thus $a \prec_2 b
\Leftrightarrow a \cup c \prec_2 b \cup c$, so $\prec_2$ is a boolean
term order. \hfill $\Box$
\end{proof}

\begin{remark} \label{sameorder} 
Note that $\prec_1$ and $\prec_2$ agree on every pair of
disjoint sets  except for the
pair $\{ \alpha, \beta \}$.
\end{remark}

\begin{example}Flipping across $\{4 \} \prec \{1,2\}$ in the term order of 
Example \ref{noncohorder} gives:
\begin{eqnarray*} 
& &\es \prec \{1 \} \prec \{2 \} \prec \{3 \} \prec \{1,2 \} \prec \{4
 \} \prec \{5 \} \prec \{1,3 \} \prec \\ & & \{2,3 \} \prec \{1,4 \}
 \prec \{1,5 \} \prec \{2,4 \} \prec \{2,5 \} \prec \{1,2,3 \} \prec
 \\ & & \{3,4 \} \prec \{1,2,4 \} \prec \{3,5 \} \prec \ldots
\end{eqnarray*}
Only the first half of the order is given, as the second half can be
deduced from the first by Corollary \ref{complement}.  This is a
coherent term order, given by the weight vector $(7,10,16,20,22)$.
\end{example}

\begin{remark} \label{fewflips} The analogue of Proposition \ref{primpairs} for flippable pairs is
false.  The following noncoherent term order in six variables, with
the flippable pairs marked by $\prec_*$, has the same order on its
flippable pairs as the coherent term order given by the vector
$(6,14,15,18,28,38 )$. It is thus not determined by the order on its
flippable pairs.

\begin{eqnarray*}
& & \es \prec \{1 \} \prec \{2 \} \prec \{1,2 \} \prec \{3 \} \prec
\{1,3 \} \prec \{4 \} \prec \{2,3 \} \prec \{1,4 \}  \prec \\ & & \{1,2,3 \} 
 \prec \{2,4 \} \prec \{5 \} \prec \{1,2,4 \} \prec \{3 ,4 \} \prec_*
\{1,5 \} \prec \{2,5 \}   \prec \\ & & \{6 \} \prec_* \{1,3,4\}  \prec
\{2,3,4 \} \prec \{1,2,5 \} \prec \{3,5 \} \prec_* \{1,6 \} \prec
\\ & & \{2,6 \} \prec \{1,2,3,4\}  \prec \{1,3,5 \} \prec \{4,5 \} \prec
\{2,3,5 \} \prec \{1,2,6 \} \prec \\ & & \{1,4,5 \} \prec_* \{3,6 \} 
\prec \{1,2,3,5 \} \prec \{1,3,6 \} \prec_* \{2,4,5\} \ldots
\end{eqnarray*}
\end{remark}

\begin{remark} The central pair is flippable in every term order.
By Corollary \ref{complement} the  $2^{n-1}+1$st term is the
complement  of the $2^{n-1}$th term, so there are no
nontrivial multiples to consider.
\end{remark}

\section{Hyperplane Arrangements and the Root System $B_n$ }

A hyperplane arrangement in $\mathbb R^n$ partitions $\mathbb R^n$
into relatively open regions of points such that in each region all
the points lie either on, or on the same side of, each hyperplane.
Let $\mathcal H_n$ denote the hyperplane arrangement consisting of all
the hyperplanes with normals in $\{ 0,1, -1 \}^n \setminus {\bf 0}^n$.
The equivalence classes of real weight vectors $w$ which determine the
same coherent generalized partial term order correspond to regions of
$\mathcal H_n$.

\begin{lemma} \label{cohgeom}
\begin{enumerate}
\item The $n$-dimensional regions of $\mathcal H_n$ are in bijection with
coherent generalized term orders.
\item Flipping across a flippable pair from one coherent term order to
another coherent term order corresponds geometrically to passing from
one $n$-dimensional region of $\mathcal H_n$ to an adjacent region.
\end{enumerate}
\end{lemma}
\begin{proof} Part 1 is immediate from the definition of $\mathcal H_n$.
Part 2 follows from Remark \ref{sameorder}, since if two regions of
$\mathcal H_n$ are on the same side of all but one hyperplane, as is
the case for the regions corresponding to term orders connected by a
flip, then the two regions must be adjacent.  

\hfill $\Box$
\end{proof}

It would thus be interesting to know the number of regions of
$\mathcal H_n$.  Unfortunately this does not appear to be a simple
combinatorial function.  One way to compute the number of regions of a
hyperplane arrangement is via its characteristic polynomial
$\chi_{\mathcal H_n}(x)$, which is defined in terms of the lattice of
intersections of the arrangement (see page 43 of \cite{orlikterao} for
details).  A result of Zaslavsky (Theorem 2.68 in \cite{orlikterao})
states that $|\chi_{\mathcal H_n}(-1)|$ is the number of regions of
the arrangement.  For a particular class of hyperplane arrangements,
known as free arrangements, the characteristic polynomial is known to
have integer roots.  In general, $\mathcal H_n$ is not a free
arrangement, as can be seen from the following table.

\vspace{4mm}
\begin{tabular}{ll}\hline
$n$ &$\chi_{\mathcal H_n}(x)$ \\ \hline
1 & $x-1$\\
2 & $(x-1)(x-3)$\\
3 & $(x-1)(x-5)(x-7)$\\
4 & $(x-1)(x-11)(x-13)(x-15)$\\
5 & $(x-1)(x-29)(x-31)(x^2-60x+971)$\\
6 & $(x-1)(x^5-363x^4+54310x^3-4182690x^2+165591769x $\\
& $-2691439347) $\\
7 & $(x-1)(x^6-1092x^6+518385x^4-136815000x^3+$ \\
&$21151739259x^2-1814252700708x+67379577529235)$\\ \hline
\end{tabular}
\vspace{4mm}

From part 2 of Lemma \ref{cohgeom} we know that the number of
facets of a particular $n$-dimensional region corresponds to the number of
coherent neighbors the corresponding term order has.  The arrangements
in dimension two and three are simplicial (all $n$-dimensional regions are simplices), but this is no longer the
case in dimension four.

\begin{example} The following term order has five
flippable pairs. As has long been known (see \cite{KPS}), and can be
seen from enumerative results (see Section \ref{numerical}), 
 all term orders are coherent in dimension four, so this term
order corresponds to a cone with five facets.  This is the maximal
number of facets of a region of $\mathcal H_4$.
\begin{eqnarray*} & & 
\es \prec \{1 \} \prec_* \{2 \} \prec_* \{3 \} \prec_* \{1,2 \} \prec \{1,3
\} \prec \{2,3 \} \\ & & \prec_* \{4 \} \prec_* \{1,2,3 \} \prec \ldots \\
\end{eqnarray*}
The comparisons marked $\prec_*$ are flippable pairs.  
\end{example}

\begin{remark} Since every $n$-dimensional region of $\mathcal H_n$ has at
least $n$ facets, coherent term orders have at least $n$ flippable
pairs.  The same is not true for noncoherent term orders.  
\end{remark}

\begin{theorem}
For $n \geq 6$ there are term orders with fewer than $n$ flippable
pairs.  These are examples of {\em flip deficiency} in the sense of
\cite{reiner}.
\end{theorem}
\begin{proof} The term order  on subsets of $\{1,2,3,4,5,6\}$ in Remark \ref{fewflips} of Section
\ref{flipsection} has only five flippable pairs.  Using this as a
base, we can construct term orders on subsets of $[n]$ with $n-1$
flippable pairs by setting $[k] \prec \{k+1\}$ for $6 \leq k \leq
n-1$. \hfill $\Box$
\end{proof}

We now develop the connection between $\mathcal H_n$ and the root system $B_n$.  
\begin{defn} The root system $B_n$ is the collection of vectors $\{ e_i : 1
\leq i \leq n \} \cup \{e_i-e_j, e_i+e_j : 1 \leq i < j \leq n \}$ in
${\mathbb R}^n$, where $e_i$ is the $i$th standard basis vector.
\end{defn}

\begin{lemma} \label{hbn}
$\mathcal H_n$ is the discriminantal arrangement of $B_n$.  In other
words, $\mathcal H_n$ is the collection of hyperplanes $\{ H : H$ is
spanned by a subset of $B_n$ \}.
\end{lemma} 

\begin{proof} Let $H$ be a hyperplane in $\mathcal H_n$, with normal vector
$v$.  Let $P=\{i : v_i>0 \}$, $N=\{i:v_i<0 \}$, and $Z=\{i:v_i=0\}$.
Then $H=\mathrm{span}(\{e_i | i \in Z \} \cup \{ e_i-e_j | i,j \in P
\} \cup \{ e_i-e_j | i,j \in N \} \cup \{e_i+e_j | i \in P, j \in N \}
)$.  For the other inclusion, let $v$ be the normal vector of a
hyperplane spanned by a subset $V$ of $B_n$.  For all $e_i-e_j \in V$
we have $v_i=v_j$, and for all $e_i+e_j \in V$ we have $v_i=-v_j$.
Also, whenever $e_i \in V$, we have $v_i=0$.  Suppose there exists
$i,j$ with $v_i,v_j \neq 0$ such that $v_i \neq \pm v_j$.  Then let
$w$ be the vector with $w_k=v_k$ when $v_k = \pm v_i$, and $w_k=0$
otherwise, and let $u$ be the vector with $u_k=v_k$ when $v_k= \pm
v_j$ and $u_k=0$ otherwise.  Then each vector in $V$ is perpendicular to  both $u$ and $w$, so lies in a codimension two subspace, and thus $V$ does not span a hyperplane.    \hfill $\Box$
\end{proof}

Let ${\mathcal M}(B_n)$ denote the oriented matroid of the vector
configuration $B_n$.  A wealth of information about oriented matroids
can be found in \cite{om}, along with many equivalent definitions.
For the purposes of the next result, we use the following chirotope
definition.

\begin{defn} Given a vector configuration $V$ of $d>n$ vectors in ${\mathbb
R}^n$, the {\em oriented matroid} $\mathcal M(V)$ is the function mapping
ordered subsets of $V$, each consisting of $n$ vectors, to
$\{+,0,-\}$, where each subset is mapped to the sign of its
determinant.  We identify the function with the sign vector of length
$n! {d \choose n}$ encoding the image of this map.
\end{defn}

We say an oriented matroid $\mathcal M (V)$ is {\em projectively unique} if
whenever $\mathcal M (W)= \mathcal M (V)$ for some vector configuration $W$ in
${\mathbb R}^d$  there is a {\em projective transformation}
(linear transformation plus scaling individual vectors) taking $W$ to
$V$.

\begin{lemma} \label{Bnunique}
${\mathcal M}(B_n)$ is projectively unique for $n \geq 3$.
\end{lemma}

\begin{proof} Suppose $\mathcal M (V)= \mathcal M (B_n)$.  We denote the
  vectors of $V$ by $v_e$, where $e$ is the vector in $B_n$ corresponding to $v_e$.

 We first apply a linear transformation of $V$ to move $v_{e_i}$ to
 $e_i$ for $1 \leq i \leq n$.   Note that for any vector $e \in B_n$ the
 determinant of the $n$ vectors $\{v_e\} \cup \{e_i : i \neq j \} $ is
$(-1)^{j-1}(v_e)_j$, while the determinant of the $n$ vectors $\{e\} \cup \{e_i : i \neq j \} $ is $(-1)^{j-1}(e)_j$.  So  $v_e$ is zero in
 exactly the same coordinates as $e$, and has the same sign in its nonzero coordinates.

 We can now do further combinations of linear transformations and
scaling of vectors to also move $v_{e_i+e_{i+1}}$ to $e_i+e_{i+1}$,
for $1 \leq i \leq n-1$. Lastly, we scale the other vectors $v_{e_i
\pm e_j}$ so they are of the form $e_i \pm ae_j$ for some $a$
depending on $i,j$ and the sign.  The lemma will follow if we can now
show that in each case $a=1$.  

We first show this for $v_{e_i+e_j}$ when $j-i$ is odd, and
$v_{e_i-e_j}$ when $j-i$ is even.  Consider $v_{e_i+e_j}=e_i+ae_j$,
with $j-i$ odd. The determinant of $\{e_l : l<i\} \cup
\{e_l+e_{l+1} :i\leq l \leq j-1 \} \cup \{v_{e_i+e_j} \} \cup \{e_l :
l>j \} $, which should be zero, is $a-1$, so we see that
$v_{e_i+e_j}=e_i+e_j$.  Considering the analogous determinant for
$v_{e_i-e_j}$, where $j-i$ is even we conclude that
$v_{e_i-e_j}=e_i-e_j$ for these $i$ and $j$.

	Now consider the determinant of $\{v_{e_1-e_2},
v_{e_2+e_3},v_{e_1+e_3} \} \cup \{ e_i : 4 \leq i \leq n \}$, which
should be zero.  If we write $v_{e_1-e_2}=e_1-ae_2$, and
$v_{e_1+e_3}=e_1+be_3$, then this determinant is $b-a$, so $b=a$.
Notice that we get the same result if we switch replace $v_{e_1-e_2}$
and $v_{e_2+e_3}$ by $v_{e_1+e_2}$ and $v_{e_2-e_3}$.  Considering
the analogous subdeterminants with similar adjacent subdeterminants,
we can conclude that there exists a single $a$ such that
$v_{e_i-e_{i+1}}=e_i-ae_{i+1}$ for all $1 \leq i \leq n-1$.  The
determinant of $\{v_{e_1-e_3},v_{e_1-e_2},v_{e_2-e_3} \} \cup \{e_i: 4
\leq i \leq n \}$, which again should be zero, is now $a^2-1$, so we see
that $a = \pm 1$.  The case $a=-1$ is ruled out by consideration of
the determinant $\{v_{e_1+e_2}, v_{e_1-e_2} \} \cup \{e_i : 3 \leq i
\leq n \}$, which is $-a-1$, and should be negative, so $a=1$, and thus $v_{e_i-e_{i+1}}=e_i-e_{i+1}$ for $1 \leq i \leq n-1$.

	From consideration, when $j-i$ is even, of the determinant of
$\{e_l : l < i \} \cup \{v_{e_i-e_{i+1}} \} \cup \{ v_{e_l +e_{l+1}} :
i+1 \leq l \leq j-1 \} \cup \{v_{e_i+e_j} \} \cup \{ e_l : l > j \}$,
which should be zero, but is $a-1$, we see that $v_{e_i+e_j}=e_i+e_j$.
Finally, the analogous determinant for $v_{e_i-e_j}$, when $j-i$ is odd,
yields $v_{e_i-e_j}=e_i-e_j$.
	
	Thus each vector $v_e\in B_n$ has been moved to the
 corresponding vector $e$, so there is a projective transformation
 moving $V$ to $B_n$, and thus we conclude that ${\mathcal M} (B_n)$
 is projectively unique. \hfill $\Box$
\end{proof}

\begin{corollary} \label{bijection}
There is a bijection between realizable one-element
extensions of ${\mathcal M}(B_n)$ and coherent generalized partial
term orders.
\end{corollary}

\begin{proof}
The implication of Lemma \ref{Bnunique} is that realizable one-element
extensions of $\mathcal M (B_n)$ correspond exactly to regions of the
discriminantal arrangement of $B_n$ for $n \geq 3$.  This is also true
by inspection for $n=2$.  Thus from Lemma \ref{hbn} we get a bijection
between realizable one-element extensions of $\mathcal M (B_n)$ and
regions of $\mathcal H_n$, and so Lemma \ref{cohgeom} gives a
bijection between realizable one-element extensions of ${\mathcal
M}(B_n)$ and coherent generalized partial term orders. \hfill $\Box$
\end{proof}

\section{ Oriented Matroids}
  
In the previous section we saw that coherent generalized partial term
orders correspond to realizable one-element extension of $\mathcal M
(B_n)$.  In this section we expand on this, showing that all generalized
partial term orders correspond to one-element extensions of $\mathcal M (B_n)$.  We assume
more familiarity with oriented matroids.

The content of Lemma \ref{hbn} was that the cocircuits of ${\mathcal
M} (B_n)$ are in bijection with the hyperplane normals in the set $\{
+1, 0, -1 \}^n \setminus {\bf 0}^n$.  We will represent these normals
by their sign vectors, which will be denoted by capital letters, such
as $X$.  The corresponding cocircuit will be denoted by
$\overline{X}$.  The passage from $X$ to $\overline{X}$ is as follows:
$\overline{X}_{e_i}=X_i$, $\overline{X}_{e_i+e_j}=X_i + X_j$, and
$\overline{X}_{e_i-e_j}=X_i - X_j$, where $+$ and $-$ applied to $\{
+,0,- \}$ evaluate to the sign of the corresponding operation on $\{ +1,0,-1
\}$.

  For example, if $n=3$, and $X=(++0)$, then
$\overline{X}=(++0+++0++)$ where the coordinates of $\overline{X}$ are
listed in the following order: the $e_i$, then the $e_i+e_j$ in
lexicographic order, and finally the $e_i-e_j$, also in lexicographic
order.

\begin{defn} $X^+$ is the set $\{i:X_i=+ \}$. $X^-$ is the set $\{i:X_i=-\}$.
Given two sign vectors $X$ and $Y$ we say that a sign vector $Z$ is an {\em
elimination candidate} for $X$ and $Y$ if $ \overline{Z}^+ \subseteq \overline{X}^+ \cup \overline{Y}^+ $ and $\overline{Z}^- \subseteq \overline{X}^- \cup \overline{Y}^-$
\label{elimcandidate}

 \end{defn}

We first describe some  $Z$ which are  elimination candidates for given $X$ and $Y$.
We can
decompose $X^{\pm},Y^{\pm}$ by writing
\begin{eqnarray*}
X^- & = & a \cup  m \cup  x \\ X^+ & = & b \cup  p \cup  y \\
Y^- & = & c \cup  m \cup  y \\ Y^+ & = & d \cup  p \cup  x \\
\end{eqnarray*}
where $a,b, c, d,m,p,x$ and $y$ are all pairwise disjoint.

\begin{lemma} \label{elims}
With decomposition as above, the following pairs $(Z^+,Z^-)$ 
are elimination candidates for $X$ and $Y$:
$$
\begin{array}{l}
(p,m) \\
 (b \cup p,a \cup m) \\
(d \cup p,c \cup m) \\
(  b \cup d \cup p \cup y,  a \cup c \cup m \cup x )\\
(  b \cup d \cup p \cup x, a \cup c \cup m \cup y )\\
\end{array}
$$

In addition, if $m=p=\es$, then
$(b \cup d,a \cup c)$ is also an elimination candidate.
\end{lemma}

\begin{proof} We can encode each of these vectors by a sign vector $S$
listing the sign of $Z_i$ for $i \in b,d,p,$ and $x$.  For example, $(p,m)$ can
be encoded as $(0,0,+,0)$.  This coding is reversible, as for all
these pairs $(Z^+,Z^-)$, and also for $(Z^+,Z^-) \in
\{(X^+,X^-),(Y^+,Y^-)\}$, we have $(b \cup d \cup p) \cap Z^- =
\emptyset$, and $A \subseteq Z^+ \Rightarrow B \subseteq Z^-$, for $(A,B) \in \{
(b,a), \, (c,d), \, (p,m)$, $(x,y), \, (y,x) \}$.  In this example,
when we pass to $\overline{(p,m)}$, this can be encoded as
$(0,0,+,0,0,+,0,+,0,+,0,-,0,-,0,+)$, with the coordinates ordered first
$b,d,p,x$, then sums in lexicographic order, and finally differences in
lexicographic order.  The elimination candidate condition is then
equivalent to requiring that, under this encoding, whenever
$\overline{S}_e=+$, either $\overline{X}_e=+$ or $\overline{Y}_e=+$,
and similarly whenever $\overline{S}_e=-$ either $\overline{X}_e=-$ or
$\overline{Y}_e=-$.

This is easy to check.  We have the following encoding of $X$ and $Y$:
\begin{eqnarray*}
\overline{X} & = & (+,0,+,-,+,+,0,+,-,0,+,0,+,-,+,+) \\
\overline{Y} & = & (0,+,+,+,+,+,+,+,+,+,-,-,-,0,0,0) \\
\end{eqnarray*}
From this we can see that $(p,m)$ is an elimination candidate.  The rest of the proof is listing the remaining four vectors.  

In the case where $m=p=\es$, we have the reduced encoding of $X$ and
$Y$ as $(+,0,-)$ and $(0,+,+)$, leaving out the $p$th coordinate.  The
condition to be an elimination candidate remains the same, so from the
encoding of $(b \cup d, a \cup c)$ as $(+,+,0,+,+,+,0,+,+)$ we can see that it is an elimination candidate. \hfill  $\Box$
\end{proof}

We are now ready for the main theorem.  A generalized
partial term order $\prec$ induces a function $\mu$ from $\{+,0,- \}^n$ to
$\{+,0,-\}$ by

$$\mu(X)=\left\{ \begin{array}{ll} + & \mathrm{ if } \, \, X^- \prec
X^+ \\ - & \mathrm{ if } \, \, X^+ \prec X^- \\ 0 & 
\mathrm{ otherwise} \\
\end{array} \right.
$$

Let $\sigma$ be the induced function on the cocircuits of ${\mathcal
M} (B_n)$ given by $\sigma(\overline{X})=\mu(X)$.

\begin{theorem} \label{toeqexts}
The map $\sigma$ is a localization, and so the set of  generalized partial term orders is in bijection with a subset of the set of one-element extensions of
${\mathcal M}(B_n)$.  \end{theorem}

\begin{proof} Since  different generalized partial term orders give different
functions $\mu$, it suffices to prove that $\sigma$ is a localization.
By Corollary 7.1.9 of \cite{om}, this is equivalent to showing that
the set $\sigma^{-1}(\{+,0 \})$ satisfies the weak cocircuit
elimination axiom.  This says that for any $X$,$Y$ with $\overline{X}
\neq -\overline{Y}$ and $\mu(X),\mu(Y) \in \{+,0 \}$ such that
$\overline{X}_e=+$ and $\overline{Y}_e=-$ for some $e \in B_n$, there
is some $Z$ with $\overline{Z}_e=0$ and $\mu(Z) \in \{+,0 \}$ which is
an elimination candidate for $X$ and $Y$.

 We first show that if $\mu(X), \mu(Y) \in
\{+,0\}$ then there is an elimination candidate, $Z$, for $X$ and $Y$
with $\mu(Z) \in \{+,0 \}$.  We decompose
$X^{\pm},Y^{\pm}$ as above, and write the candidate $Z$ as
$(Z^+,Z^-)$. Parentheses will be omitted at times to simplify
notation.  Note that the condition that $\mu(X),\mu(Y) \in \{+,0\}$
means $a \cup m \cup x \preceq b \cup p \cup y$ and $ c \cup m \cup y
\preceq d \cup p \cup x$, and also that $\overline{X} \neq - \overline{Y}$.  The argument divides into two cases.  We
make repeated use of Lemma \ref{multlemma}.

\begin{itemize}
\item[Case I:] $m \preceq p$.  Then $\mu(p,m) \in \{+,0\}$.
\item[Case II:] $p \prec m$ or $m=p=\emptyset$.  Then $a \cup x \preceq
 	b \cup y$, and $c \cup y \preceq d \cup x$.  There are three
 	further cases.  
	
	\begin{itemize}
 	\item[Case (a):] $x \prec y$.  Then $c \cup m \prec
d \cup p$, so $\mu(d \cup p,c
\cup m)=+$.  Also $c \prec d$, so
$\mu(b \cup d \cup p \cup y,a
\cup c \cup m \cup x)=+$.  Note that these are both nonempty pairs in the case $m=p=\emptyset$.

 	\item[Case (b):] $y \prec x$.  Then
$\mu(b \cup p,a \cup m)=+$ and
$a \prec b$, so \\ $\mu( b
\cup d \cup p \cup x,a \cup c \cup m \cup
y)=+$. Again, these are both nonempty pairs.

	\item[Case (c):] $x \sim y$ or $x=y=\es$.  We break into two
		further cases.  ]`
		\begin{itemize} 

		\item[Case (i):] $p \prec m$. Then $a \cup m \preceq b
		\cup p$, so $\mu( b \cup p,a \cup m) \in \{+,0 \}$,
		and similarly $\mu(d \cup p,c \cup m) \in \{+,0 \}$.
		Also $c \prec d$, so $\mu(b \cup d \cup p \cup y,a
		\cup c \cup m \cup x)=+$.  
		
		\item[Case (ii):] $m=p= \es$. Then $\mu(b \cup d,a
		\cup c)$, $\mu ( b \cup d \cup y,a \cup c \cup x)$, $
		\mu(b \cup d \cup x,a \cup c \cup y) \in
		\{+,0\}$. Since $\overline{X} \neq -\overline{Y}$,
		these are all nonempty pairs.  
		\end{itemize}

\end{itemize}

\end{itemize}

\begin{figure}

\psfrag{a}{\footnotesize{$m \preceq p$}}
\psfrag{b}{\tiny{(1) $(p,m)$}} 
\psfrag{c}{\footnotesize{$p \prec m$ or }}
\psfrag{d}{\footnotesize{    $m=p=\emptyset$}}
\psfrag{e}{\footnotesize{$x \prec y$}}
\psfrag{f}{\tiny{(1) $(d \cup p, c \cup m)$}}
\psfrag{g}{\tiny{(2) $(b \cup d \cup p \cup y,$}}
\psfrag{h}{\tiny{ \, \, \, \,  $ a \cup c \cup m \cup x)$}} 
\psfrag{i}{\footnotesize{$y \prec x$}}
\psfrag{j}{\tiny{(1) $(b \cup p, a \cup m)$}}
\psfrag{k}{\tiny{(2) $( b \cup d \cup p \cup x,$}}
\psfrag{l}{ \tiny{\, \, \, \, $ a \cup c \cup m \cup  y)$}}
\psfrag{m}{\footnotesize{$y \sim x$ or }}
\psfrag{mm}{\footnotesize{$x=y=\emptyset$}}
\psfrag{n}{\tiny{(1) $(b \cup p, a \cup m)$}}
\psfrag{o}{\tiny{(2) $(d \cup p, c \cup m)$}}
\psfrag{p}{\tiny{(3) $( b \cup d \cup p \cup y,$}}
\psfrag{q}{ \tiny{\, \, \,$ a \cup c \cup m \cup x)$}}
\psfrag{r}{\footnotesize{$x \prec y$}}
\psfrag{s}{\tiny{(1) $(d,c)$}}
\psfrag{t}{\tiny{ (2) $(b \cup d \cup y,$}}
\psfrag{u}{\tiny{ \, \, \, \, $a \cup c \cup x)$}}
\psfrag{v}{\footnotesize{$y \prec x$}}
\psfrag{w}{\tiny{(1) $(b,a)$}}
\psfrag{x}{\tiny{(2) $(b \cup d \cup x,$}}
\psfrag{y}{\tiny{ \, \, \, \, $a \cup c \cup y)$}}
\psfrag{z}{\footnotesize{$y \sim x$ or $x=y=\emptyset$}}
\psfrag{aa}{\tiny{(1) $(b \cup d, a \cup c)$}}
\psfrag{bb}{\tiny{(2) $(b \cup d \cup y,$}}
\psfrag{cc}{\tiny{ \, \, \, \, $a \cup c \cup x)$}}
\psfrag{dd}{\tiny{(3) $(b \cup d \cup x,$}}
\psfrag{ee}{\tiny{ \, \, \, \, $a \cup c \cup y)$}}
\psfrag{ff}{\footnotesize{$p \prec m$}}
\psfrag{gg}{\footnotesize{$m=p=\emptyset$}}

\includegraphics{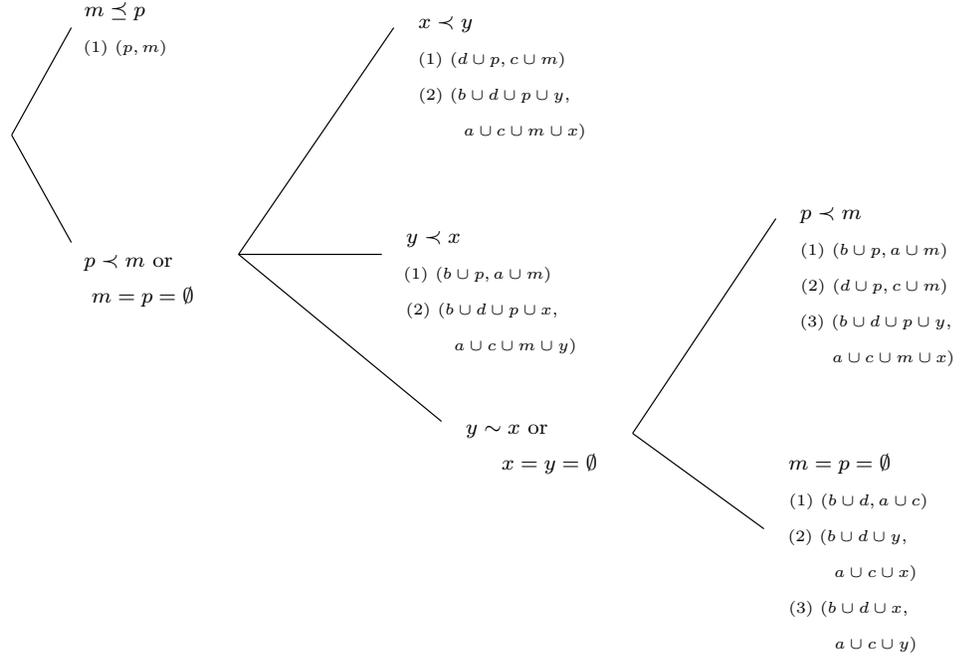}
\caption{$\mu$-nonnegative vectors}
\end{figure}

So we have the diagram of cases shown in Figure 1.  Note that in each
branch the pair is nonempty, so it does represent a cocircuit.  By
Lemma \ref{elims} the pair in each branch is an elimination candidate
for $X$ and $Y$. So we have shown that if $\mu(X),\mu(Y) \in \{+,0\}$
then there is an elimination candidate, $Z$, for $X$ and $Y$ with
$\mu(Z) \in \{+,0\}$.  We are now ready to show that
$\sigma^{-1}(\{+,0\})$ satisfies the weak cocircuit elimination
condition.  Recall that this involves showing for all $X,Y \in
\mu^{-1}(\{+,0\})$ with $\overline{X} \neq -\overline{Y}$, and $e \in
B_n$ such that $\overline{X}_e=+=-\overline{Y}_e$ there exists a $Z
\in \mu^{-1} (\{+,0\})$ which is an elimination candidate for $X$ and
$Y$ with $\overline{Z}_e=0$.  To do this, it suffices to show that if
$\overline{X}_e=+$ and $\overline{Y}_e=-$ then on every branch of
Figure 1 there is a $Z$ with $\overline{Z}_e=0$.  This is again a
consideration of cases, depending on the form of $e$. The following
tables enumerate the cases, with the rightmost column indicating which
$Z$ is chosen for each branch of the above diagram.  We read the nodes
of the diagram in the order the corresponding cases appear above; that
is, top to bottom then left to right.  The numbers are those in the
diagram, so $11111$ represents the choice of $(p,m)$, $(d \cup p, c
\cup m)$, $(b \cup p, a \cup m)$, $(b \cup p, a \cup m)$, and $(b \cup
d, a \cup c)$ in the corresponding branches.

\begin{enumerate}

\item $e=e_i$.  Then $i \in x \cup
y$ and we require $Z_i=0$, so we take $11111$.

\item $e=e_i+e_j$.  Then we require $Z_i=-Z_j$

\begin{tabular}{rrrrccl} \hline

$X_i$ & $ X_j $ & $ Y_i $ & $ Y_j $ & $ i \in $ & $ j \in $ & $  Z$ \\ \hline
$+ $ & $ + $ & $ - $ & $ - $ & $ y $ & $ y $ & $11111  $  \\
$+ $ & $ + $ & $ - $ & $ 0 $ & $ y $ & $ b $ & $11223$  \\
$+ $ & $ + $ & $ 0 $ & $ - $ & $ b $ & $ y $ & $11223$  \\ 
$+ $ & $ 0 $ & $ - $ & $ - $ & $ y $ & $ c $ & $12112$    \\
$+ $ & $ 0 $ & $ - $ & $ 0 $ & $ y $ & $ 0 $ &  $11111  $   \\
$+ $ & $ 0 $ & $ 0 $ & $ - $ & $ b $ & $ c $ & $12231$\\
$0 $ & $ + $ & $ - $ & $ - $ & $ c $ & $ y $ & $12112$  \\
$0 $ & $ + $ & $ - $ & $ 0 $ & $ c $ & $ b $ &$12231$\\
$0 $ & $ + $ & $ 0 $ & $ - $ & $ 0 $ & $ y $ & $11111  $   \\ \hline

\end{tabular}

\item $e=e_i-e_j$.  Then we require $Z_i=Z_j$.

\begin{tabular}{rrrrccl} \hline

$X_i$ & $ X_j $ & $ Y_i $ & $ Y_j $ & $ i \in  $ & $ j \in  $ & $  Z$ \\ \hline
$+ $ & $ - $ & $ - $ & $ + $ & $ y $ & $ x $ & $11111  $ \\
$+ $ & $ - $ & $ - $ & $ 0 $ & $ y $ & $ a $ & $11223$     \\
$+ $ & $ - $ & $ 0 $ & $ + $ & $ b $ & $ x $ & $11223$      \\ 
$+ $ & $ 0 $ & $ - $ & $ + $ & $ y $ & $ d $ & $12112$  \\
$+ $ & $ 0 $ & $ - $ & $ 0 $ & $ y $ & $ 0 $ & $11111  $    \\
$+ $ & $ 0 $ & $ 0 $ & $ + $ & $ b $ & $ d $ & $12231$\\
$0 $ & $ - $ & $ - $ & $ + $ & $ c $ & $ x $ & $12112$  \\
$0 $ & $ - $ & $ - $ & $ 0 $ & $ c $ & $ a $ & $12231$ \\
$0 $ & $ - $ & $ 0 $ & $ + $ & $ 0 $ & $ x $ & $11111  $   \\ \hline

\end{tabular}

\end{enumerate}

This shows that $\sigma^{-1}(\{+,0\})$ satisfies the weak cocircuit
elimination axiom, so $\sigma$ is a localization.  \hfill $\Box$
\end{proof}

\begin{corollary}
A noncoherent generalized partial term order determines a nonrealizable one-element extension of $\mathcal M (B_n)$.
\end{corollary}

\begin{proof} This follows directly from Theorem \ref{toeqexts} and Corollary
\ref{bijection}. \hfill $\Box$.
\end{proof}

 We will show in Remark \ref{noconverse} that there are one-element
extensions of $\mathcal M (B_n)$ which are not induced by a
generalized partial term order in this way.  To this end we
characterize those one-element extensions which are in the image of
the bijection.

\begin{proposition} Let $\mu$ be a function from $\{ +,-,0 \}^n \setminus {\bf 0}^n$ to
$\{+,0,- \}$.  It is induced by a generalized partial term order if
and only if it satisfies the following criteria:
\begin{enumerate}
\item $\mu(-x)=-\mu(x)$ 
  \item {\bf (First Addition Condition)} If $\mu(x)=\mu(y)=+$, and $x_i
\neq y_i $ whenever $x_i \neq 0$ for $1 \leq i \leq n$, then $\mu(z )=+$, where
$$z_i=\left\{ 
\begin{array}{lll} x_i & y_i=0 \\ y_i & x_i=0
\\ 0 & x_i=-y_i \neq 0 \\ \end{array} \right.$$ 

\item {\bf (Second Addition Condition)} If $\mu(x)=0$, $\mu(y)=+$,
and $x_i \neq y_i$ whenever $x_i \neq 0$ for $1 \leq i \leq n$, then $\mu(z )=+$, where $z$ is as above.
\end{enumerate}
\end{proposition}

\begin{proof} The necessity of the first condition is immediate from the way
$\mu$ is induced by a term order. Lemma \ref{multlemma} implies the
second and third conditions.

For sufficiency, let $\mu$ be a function satisfying the above
conditions.  Given two subsets of $[n]$, define their order to be the one
induced on the pair of subsets obtained by removing their intersection from each subset.  Suppose that this does not define a partial
order on the set of all subsets of $[n]$.  Then there is some string of
inequalities $b_0 \prec b_1 \prec b_2 \prec \ldots \prec
b_n \preceq b_0$.  We can reduce this to a string of
three distinct monomials.  If $b_{i} \prec b_{i+2}$, we can remove $b_{i+1}$ from the chain and repeat this procedure with the shorter chain.  Otherwise we have the three-element chain  $b_i \prec b_{i+1} \prec b_{i+2} \preceq
b_i$.
Denote this three-element chain $\alpha \prec \beta \prec
\gamma \preceq \alpha$. We may assume that $\alpha \cap
\beta \cap  \gamma = \emptyset$.  Write 
\begin{eqnarray*}
\alpha & = & \alpha^{\prime} \cup x \cup y \\ \beta & = &
\beta^{\prime} \cup x \cup z \\ \gamma & = & \gamma^{\prime} \cup y
\cup z\\ \end{eqnarray*} where $\alpha^{\prime}$, $\beta^{\prime}$,
$\gamma^{\prime}$, $x$, and $y$ are pairwise disjoint.  Then we have
$\alpha^{\prime} \cup y \prec \beta^{\prime} \cup z$, $\beta^{\prime}
\cup x \prec \gamma^{\prime} \cup y$, and $\gamma^{\prime} \cup z
\preceq \alpha^{\prime} \cup x$.  But the addition condition demands
$\alpha^{\prime} \cup x \prec \gamma^{\prime} \cup z$ so $\mu$
determines a partial order on subsets of $[n]$.

Suppose the partial order does not satisfy the condition on
incomparable elements, so there is $\alpha \sim \beta$, $\alpha \prec
\gamma$, but $\gamma \preceq \beta$.  Decompose $\alpha, \beta$, and
$\gamma$ as above.  Then $\alpha^{\prime} \cup y \sim \beta^{\prime}
\cup z$, and $\alpha^{\prime} \cup x \prec \gamma^{\prime} \cup z$ but
$\gamma^{\prime} \cup y \preceq \beta^{\prime} \cup x$.  Then
$\mu(\alpha^{\prime} \cup y,\beta^{\prime} \cup z)=0$,
$\mu(\gamma^{\prime} \cup z,\alpha^{\prime} \cup x)=+$, and
$\mu(\gamma^{\prime} \cup y,\beta^{\prime} \cup x) \in \{0,- \}$.
However, the second addition condition demands $\mu(\gamma^{\prime}
\cup y,\beta^{\prime} \cup x)=+$, so the order satisfies the condition
on incomparable elements.

	 Since the partial order satisfies the first condition of
Definition \ref{gpto} (the multiplicative condition) by construction,
it is a generalized partial term order, and so $\mu$ is induced by a
generalized partial term order.  \hfill $\Box$
\end{proof}

\begin{remark} \label{noconverse}  The following list of hyperplanes in 
${\mathbb R}^3$ satisfies the weak cocircuit elimination axiom when considered 
as cocircuits of $B_3$: 
$$\begin{array}{l}(-+0),(0-+),(+0-),(+00),(0+0),(00+),(0++),\\
(+0+),(++0),(+++),(++-),(+-+),(-++
)\\
\end{array}$$

Thus when we set $\sigma(\overline{X})=+$ for each hyperplane $X$,
this determines a one-element extension of $\mathcal M (B_3)$.  There
is no term order, however, which induces this $\sigma$, as the set
does not satisfy the first addition condition.  Explicitly, if these
were the positive cocircuits from some term order, then from the first
three we have $x \prec y$, $y \prec z$, and $z \prec x$, which is a
contradiction.  So not all one-element extensions of $\mathcal M
(B_n)$ are induced by generalized partial term orders.
\end{remark}

\section{Numerical Results and Examples} \label{numerical}

 In \cite{FineGill} Fine and Gill give crude bounds on the number of
antisymmetric comparative probability relations, and the number of
antisymmetric additive comparative probability relations, and the
first terms of each sequence.  

With improved computer speeds it is now possible to evaluate a few
more terms in each sequence.  The values calculated are displayed in
the following table.   The numbers in the second and third
column are divided by $n!$, taking into account the action of the
symmetric group.

\vspace{4mm}
\begin{tabular}{rrr} \hline
	$n$ & Number of term orders$/n!$  & Number of coherent term orders$/n!$\\ \hline
	1 & 1 &1\\
	2 & 1 &1\\
	3& 2 & 2 \\
	4 & 14 & 14\\
	5 & 546 & 516 \\
	6 & 169444 & 124187\\
	7 & 560043206 & 214580603 \\ \hline
\end{tabular}
\vspace{4mm}

These were calculated with the aid of programs written by Michael
Kleber and Josh Levenberg.  To calculate the total number of term
orders we used a recursive procedure, calculating for each term order
on subsets of $[n-1]$ the number of ways the subsets involving $n$
could be shuffled in.  To calculate the number of coherent term orders
we enumerated the regions of the corresponding hyperplane arrangement
$\mathcal H_n$.  These numbers have since been checked directly for $n
\leq 6$.

The following table enumerates the number of flippable pairs each term
order has for orders on subsets of $\{1,2,3,4,5,6\}$.

\vspace{4mm}
\begin{tabular}{rrr} \hline
$k$ & Number with $k$ flippable pairs \\ \hline
5 & 107 \\
6 & 14699 \\
7 & 46626\\
8 & 56707\\
9 & 35555\\
10 & 12763\\
11 & 2633\\
12 & 334\\
13 & 20\\ \hline
Total & 169444 \\ \hline
\end{tabular}
\vspace{4mm}

Notice that there are some orders with significantly more than six
flippable pairs.

\begin{defn} The {\em Baues poset} for term orders has as its elements 
all generalized partial term orders.  The order relation is that a
term order $\prec_1$ is less than another term order $\prec_2$ if
$\prec_1$ is a refinement of $\prec_2$.  This poset has a $\hat{1}$, the generalized partial term order with all subsets unrelated.
\end{defn}

\begin{proposition}
There exists a noncoherent boolean term order which lies below no coherent partial term order in the Baues poset except $\hat{1}$.
\end{proposition}
\begin{proof} The following order is such an example:
\begin{eqnarray*}
& & \es \prec \{1 \} \prec \{2 \} \prec \{1,2 \} \prec \{3 \} \prec
\{1,3 \} \prec \{2,3 \} \prec \{1,2,3 \} \prec \\ & & \{4 \} \prec
\{1,4 \} \prec \{2,4 \} \prec \{1,2,4 \} \prec \{3,4 \} \prec \{5 \}
\prec \{1,3,4 \} \prec \\ & & \{2,3,4 \} \prec \{1,5 \} \prec \{2,5 \}
\prec \{1,2,3,4 \} \prec \{1,2,5 \} \prec \{3,5 \} \prec \\ & & \{1,3,5 \}
\prec  \{2,3,5 \} \prec \{6 \} \prec \{1,2,3,5 \} \prec \{1,6 \}
\prec \{4,5 \} \prec \\ & &  \{1,4,5 \} \prec \{2,6 \} \prec  \{1,2,6 \}
\prec \{3,6 \} \prec \{1,3,6 \} \ldots \\
\end{eqnarray*}
If it were less than a coherent partial term order other than $\hat{1}$ in
 the Baues poset, the linear program obtained by seeking a weight
 vector with all $\prec$ above relaxed to $\preceq$ would have a nonzero solution.  This is not the case.  \hfill $\Box$
\end{proof}

Note that the example in the proof above is adjacent to a coherent
order, however.  Flipping over $\{5 \} \prec \{1,3,4\}$ yields the
coherent order given by the vector $(2,9,12,28,48,70)$.

The following term order in six variables, with flippable pairs marked
with $\prec_*$, is not adjacent to any coherent term order:

\begin{eqnarray*}
& & \es \prec_* \{1 \} \prec \{2 \} \prec \{1,2 \} \prec \{3 \}
\prec \{1,3 \} \prec \{4 \} \prec \{1,4 \} \prec \{5 \}  \prec
\\ & & \{1,5 \} \prec* \{2,3 \} \prec \{1,2,3 \} \prec \{6 \} \prec \{1,6 \}
\prec \{2,4 \} \prec \{1,2,4 \} \prec \\ & & \{3,4 \} \prec \{1,3,4 \}
\prec_* \{2,5 \} \prec \{1,2,5 \} \prec \{2,6 \} \prec  \{1,2,6
\}  \prec_* \\ & & \{3,5 \} \prec \{1,3,5 \} \prec \{3,6 \} \prec \{1,3,6 \}
\prec_* \{4,5 \} \prec \{1,4,5 \} \prec \\ & & \{2,3,4 \} \prec \{1,2,3,4 \}
 \prec \{2,3,5 \} \prec \{1,2,3,5 \} \prec_* \{4,6 \} \ldots \\
\end{eqnarray*}

This is the case because $\{5 \} \prec \{2,3\}$, $\{3,4\} \prec
\{2,5\}$, $\{2,6 \} \prec \{3,5\}$, and $\{2,3,5\} \prec \{4,6\}$,
which, in monomial notation, gives $x_2^2x_3^2x_4x_5^2x_6 \prec
x_2^2x_3^2x_4x_5^2x_6$, implying that $\prec$ is noncoherent.
Flipping across any of the flippable pairs does not change any of
these four inequalities, so none of the neighbors of $\prec$ are
coherent.

This term order was constructed by taking a ``lexicographic product''
of the noncoherent term order of Example \ref{noncohorder} with $\es
\prec \{1\}$.  Using this product construction we can construct
boolean term orders in $n$ variables which are at least $2^{n-5}$
flips from a coherent boolean term order.

\section{Questions}
The following questions are natural in the context of a Baues problem.

\begin{itemize}

\item Is the space of term orders connected by flips?  This has been
experimentally verified for $n \leq 6$.

\item What is the homotopy type of the poset of (generalized) partial
term orders?  In particular, is it spherical?  The subposet of
coherent generalized partial term orders is easily seen to be spherical.

\item What is the limit of the ratio of the number of coherent
term orders to the total number of term orders as $n$ increases?  Is
it zero?

\item What is an upper bound for the number of coherent neighbors for
a coherent term order?  In other words, how many facets do regions of
the hyperplane arrangement have?

\item Can we give a lower bound for the total number of flippable pairs
in all term orders?  A good upper bound?

\end{itemize}

\end{document}